\documentclass[notitlepage,11pt]{article}
\usepackage{amssymb,amsmath,comment,cases}
\catcode`\@=11
\@addtoreset{equation}{section}

\catcode`\@=12

\newtheorem{Theorem}{Theorem}[section]

\newtheorem{Corollary}[Theorem]{Corollary}
\newtheorem{Remark}[Theorem]{Remark}

\def\QED{\hfill {$\square$}\goodbreak \medskip}

\def\R{{\mathbb R}}
\def\H{H}
\def\abs{{a_{*}}}

\begin{document}

\title{Symmetry breaking of extremals for the Caffarelli-Kohn-Nirenberg inequalities \\in a non-Hilbertian setting}

\author{
Paolo Caldiroli\footnote{Dipartimento di Matematica, Universit\`a di Torino, via Carlo Alberto, 10 -- 10123 Torino, Italy. Email: {\tt paolo.caldiroli@unito.it}. Partially supported by the PRIN2009 grant "Critical Point Theory and Perturbative Methods for Nonlinear Differential Equations"}
~and~
Roberta Musina\footnote{Dipartimento di Matematica ed Informatica, Universit\`a di Udine, via delle Scienze, 206 -- 33100 Udine, Italy. 
Email: {\tt roberta.musina@uniud.it}. Partially supported by Miur-PRIN 2009WRJ3W7-001 ``Fenomeni di concentrazione e {pro\-ble\-mi} di analisi geometrica''}}

\date{}

\maketitle

\begin{abstract}
We provide an explicit necessary condition to have that no extremal for the best constant in the Caffarelli-Kohn-Nirenberg inequality is radially symmetric.
\footnotesize
\medskip

\medskip

\end{abstract}
\medskip

\normalsize
\section{Introduction}
\label{S:introduction}
This paper deals with the celebrated Caffarelli-Kohn-Nirenberg inequalities in $\R^{n}$, $n\ge 2$. 
Let $p,q$ be given exponents, such that
$$
1<p<q<p^{*} ~, 
$$
where $p^{*}=\frac{np}{n-p}$ if $p<n$, and $p^{*}=\infty$ if $p\ge n$. For any
$a>p-n$ define
\begin{equation}
\label{eq:b}
 b_a=n-q\frac{n-p+a}{p}.
 \end{equation}
In \cite{CKN} and \cite{Maz} it is proved that there exists a constant $c=c(N,p,a,q)>0$ such that 
\begin{equation}
\label{eq:CKN}
c\left(\int_{\R^{n}}|x|^{-b_a}|u|^{q}~\!dx\right)^{p/q}\le\int_{\R^{n}}|x|^{a}|\nabla u|^p~\!dx
\end{equation}
for any $u\in C^\infty_c(\R^{n})$. Let 
$$
{\mathcal D}^{1,p}(\R^{n};|x|^a dx)=\left\{u\in L^p(\R^n;|x|^{a-p} dx)~|~\int_{\R^n}|x|^a|\nabla u|^p~\!dx<\infty~\right\}.
$$
Notice that $\mathcal D^{1,p}(\R^{n}; dx)=\mathcal D^{1,p}(\R^{n})$
if $p<n$ and $a=0$. The best constant in (\ref{eq:CKN}) is given by
\begin{equation}
\label{eq:minimization}
S_{p}(a,q) :=\inf_{\scriptstyle u\in{\mathcal D}^{1,p}(\R^{n};|x|^a dx)\atop\scriptstyle u\ne 0}\frac{\displaystyle\int_{\R^{n}}|x|^a|\nabla u|^p~\!dx}{\left(\displaystyle\int_{\R^{n}} |x|^{-b_a}|u|^{q}~\!dx\right)^{p/q}}.~
\end{equation}

We point out that the minimization problem (\ref{eq:minimization}) is characterized by a lack of compactness, since the integrals in (\ref{eq:CKN}) are both invariant with respect to the group of 
transforms
\begin{equation}
\label{eq:dilations}
u(x)\mapsto \rho^{H}~\!u(\rho x)\quad(\rho>0)\quad\text{where}\quad H=\frac{n-p+a}{p}~\!.
\end{equation}
Nevertheless, the infimum $S_{p}(a,q)$ is always achieved: we quote \cite{CaWa} for $p=2$ (see also \cite{CaMu99} for a very first result in a special case), 
and \cite{Ho}, \cite{GaMu1} in the  generality of the 
above assumptions.

The exponent $H$ in (\ref{eq:dilations}) has a special meaning, as it is related to the Hardy inequality:
it is well known that
$$
H^{p}=\inf_{\scriptstyle u\in{\mathcal D}^{1,p}(\R^{n};|x|^a dx)\atop\scriptstyle u\ne 0}\frac{\displaystyle\int_{\R^{n}}|x|^a|\nabla u|^p~\!dx}{\displaystyle\int_{\R^{n}} |x|^{a-p}|u|^{p}~\!dx}~\!.
$$

The presence of weights in the integrals in (\ref{eq:minimization})  affects the property of minimizers to be radially symmetric or not. In this paper we prove the following result.

\begin{Theorem}
\label{T:BS1}
Assume $p<q<p^*$ and let $a>p-n$ be given. 
If 
\begin{equation}
\label{eq:BS}
\frac{1}{n-1}\left(\frac{n-p+a}{p}\right)^2>\frac{1}{q-p}-\frac{1}{q+p'}~,
\end{equation}
where $p'$ is the conjugate exponent to $p$, then no minimizer for $S_{p}(a,q)$ is radial.
\end{Theorem}

As a consequence to Theorem \ref{T:BS1}, we immediately get a multiplicity result for problem
\begin{equation}
\label{eq:equation}
\begin{cases}
-\mathrm{div}(|x|^a|\nabla u|^{p-2}\nabla u)=|x|^{-b_a}|u|^{q-2}u\text{ ~on~ $\R^{n}$}\\
\displaystyle\int_{\R^{n}}|x|^a|\nabla u|^p~\!dx<\infty.
\end{cases}
\end{equation}
\begin{Corollary}
If the assumptions in Theorem \ref{T:BS1} are satisfied, then 
problem (\ref{eq:equation}) has at least two distinct nonnegative and nontrivial solutions.
\end{Corollary}

Theorem \ref{T:BS1} implies in particular that for any fixed $p,q$, if 
$a$ is larger than the unique $a^*>p-n$ such that
$$
\frac{1}{n-1}\left(\frac{n-p+a^*}{p}\right)^2=\frac{1}{q-p}-\frac{1}{q+p'}~,
$$
then symmetry breaking occurs.  The next result better explains this phenomenon from a qualitative point of view (see also Remark \ref{R:np} for the case
$n\le p$).

\begin{Theorem}
\label{T:2}
Let $p,q$ be given, with $n>p$ and $1<p<q<p^*$.
There exists an exponent $\abs\ge 0$ such that 
\begin{itemize}
\item[(i)]
for any $a>\abs$ no extremal for $S_{p}(a,q)$ is radially symmetric;
\item[(ii)]
for any $a\le \abs$ there exists a radial extremal for $S_{p}(a,q)$.
\end{itemize}
\end{Theorem}
By Theorem \ref{T:BS1} we have that $\abs\le a^*$. We conjecture that condition (\ref{eq:BS}) can not be improved, that is, $\abs=a^*$.

More results are available in the Hibertian case  $p=2$. 
Symmetry breaking was already observed by Catrina and Wang in \cite{CaWa}. Then in \cite{FeSc} Felli and Schneider gave a sharper description of the region in which symmetry breaking occurs. More precisely, from Corollary 1.2 of \cite{FeSc} one gets that if 
$$
\frac{1}{n-1}\left(\frac{n-2+a}{2}\right)^2>\frac{1}{q-2}-\frac{1}{q+2}
$$
then no minimizer for $S_{2}(a,q)$ is radial. Notice that condition (\ref{eq:BS}) with $p=2$ coincides with the Felli-Schneider's one. 

As concerns symmetry results, 
in the recent paper \cite{DEL} Dolbeault, Esteban and Loss prove that minimizers for $S_{2}(a,q)$ are necessarily radial when
$$
\frac{1}{n-1}\left(\frac{n-2+a}{2}\right)^2\le\frac{1}{q-2}-\frac{1}{4}.
$$
In particular, $a^*>0$ in the Hilbertian case $p=2$.

In \cite{BW} and \cite{SW} the case $p>1$ is considered. The main result in those papers gives symmetry breaking for any $a\ge a_0$, with $a_0=a_0(n,p,q)$ large enough. However, no explicit estimate for $a_0$ is given there.  
The proofs in \cite{BW} and in \cite{SW} are based on the analysis of the asymptotic behaviour of the best constants $S_p(a,q)$ and $S_{p,{\rm rad}}(a,q)$ (on radial functions) as $a\to\infty$.

Our proof relies on the explicit knowledge of radial solutions to (\ref{eq:equation}) and is developed by computing the second variation of the functional defined by the quotient in (\ref{eq:minimization}), along a suitably chosen direction $v\in {\mathcal D}^{1,p}(\R^{n};|x|^a dx)$.

\section{Proof of Theorem \ref{T:BS1}}
It has been recently proved in \cite{MuAMPA} that, up to change of sign and up to the transform (\ref{eq:dilations}), problem (\ref{eq:equation}) has a unique nontrivial radial solution $U$ in ${\mathcal D}^{1,p}(\R^{n};|x|^a dx)$. More precisely, $U$ is given by
\begin{equation}
\label{eq:U_explicit}
U(x)= C\left(1+|x|^{\frac{(n-p+a)(q-p)}{p(p-1)}}\right)^{-\frac{p}{q-p}}\quad\text{where}\quad C=\left(\frac{q}{p}~\!\frac{\left({n-p+a}\right)^p}{(p-1)^{p-1}}\right)^{\frac{1}{q-p}}.
\end{equation}
Since the function $U$ in (\ref{eq:U_explicit}) solves (\ref{eq:equation}), then $U$ is a critical point of the functional
$$
\mathcal R(u)=\frac{\displaystyle\int_{\R^{n}}|x|^a|\nabla u|^p~\!dx}
{\displaystyle\left(\int_{\R^{n}}|x|^{-b_a}|u|^q~\!dx\right)^{p/q}}~,\quad~~\mathcal R
\colon{\mathcal D}^{1,p}(\R^{n};|x|^a dx)\setminus\{0\}\to \R~\!.
$$
Theorem \ref{T:BS1} follows by analyzing the linearized operator around
the radial extremal $U$, to  show that  if $U$ achieves the best constant $S_{p}(a,q)$
then (\ref{eq:BS}) can not hold.

Assume that $U$ is a local minimum for
$\mathcal R$. Let $\varphi_1\in H^1(\mathbb{S}^{n-1})$ be the eigenfunction of the Laplace operator on $\mathbb{S}^{n-1}$ relatively to the smaller positive eigenvalue, that is,
\begin{equation}
\label{eq:f1}
\int_{\mathbb{S}^{n-1}} \varphi_1~\!d\sigma=0~,\quad
\int_{\mathbb{S}^{n-1}}|\varphi_1|^2~\!d\sigma=1~,\quad \int_{\mathbb{S}^{n-1}}|\nabla_\sigma \varphi_1|^2~\!d\sigma=n-1~.
\end{equation}
Fix any nontrivial and radial function $v\in C^\infty_c(\R^{n}\setminus\{0\})$. Clearly, one has that $v\varphi_1\in{\mathcal D}^{1,p}(\R^{n};|x|^a dx)$. Since $U$ achieves the best constant $S_{p}(a,q)$, then the real function
${\mathcal R}_{v}\colon t\mapsto~\mathcal R(U+tv\varphi_1)$
has a local minimum at $t=0$. Notice that ${\mathcal R}_{v}$ is twice differentiable at $0$, as $\nabla U$ never vanishes on the support of $v$ (some care is needed if $p<2$). Using
$$
|\nabla(v\varphi_1)|^2=|\nabla v|^2 \varphi_1^2+|x|^{-2}|v|^2|\nabla_\sigma\varphi_1|^2,
$$
(\ref{eq:f1}) and  ${\mathcal R}_{v}''(0)\ge 0$  one gets
\begin{equation}
\label{eq:Q''}
\begin{split}
(p-1)\int_{\R^{n}}|x|^a|\nabla U|^{p-2}|\nabla v|^2~\!dx+
(n-1)&\int_{\R^{n}}|x|^{a-2}|\nabla U|^{p-2}|v|^2~\!dx\\
&\ge(q-1) \int_{\R^{n}} |x|^{-b_a}U^{q-2}|v|^2~\!dx~\!.
\end{split}
\end{equation}
Our next goal is to take
\begin{equation}
\label{eq:v}
v=|x|^{\beta H}~U^\frac{q}{p}\quad\text{where}\quad
\H=\frac{n-p+a}{p}
\end{equation}
and $\beta$ to choose in a suitable way. This step can be done by a standard limit argument provided that the exponent $\beta$ is such that the mappings $|x|^{a}|\nabla U|^{p-2}|\nabla v|^2$, $|x|^{a-2}|\nabla U|^{p-2}|v|^2$, and $|x|^{-b_a}U^{q-2}|v|^2$, with $v$ given by (\ref{eq:v}), are integrable on $\R^{n}$. To this purpose we set
\begin{equation*}
\begin{split}
&I_{0}=\int_{\R^{n}}|x|^{a+2\beta H}U^{\frac{2q}{p}-2}|\nabla U|^{p}~\!dx\\
&I_{1}=\int_{\R^{n}}|x|^{a+2\beta H-1}U^{\frac{2q}{p}-1}|\nabla U|^{p-1}~\!dx\\
&I_{2}=\int_{\R^{n}}|x|^{a+2\beta H-2}U^{\frac{2q}{p}}|\nabla U|^{p-2}~\!dx
\end{split}
\end{equation*}
and we notice that the summability at the origin and at infinity of the functions contained in the above integrals is ensured provided that $\beta$ satisfies
\begin{equation}
\label{eq:summability}
\begin{cases}
2\beta+p+(p-2)Q>0\\
2\beta<\left(2+p+\dfrac{pq}{q-p}\right)Q-p
\end{cases}
\end{equation}
respectively. 
In this case one has that
\begin{equation*}
\begin{split}
&\int_{\R^{n}}|x|^a|\nabla U|^{p-2}|\nabla v|^2~\!dx=\left(\frac{q}{p}\right)^{2}I_{0}-\frac{2\beta H q}{p}I_{1}+(\beta H)^{2}I_{2}\\
&\int_{\R^{n}}|x|^{a-2}|\nabla U|^{p-2}|v|^2~\!dx=I_{2}\\
&\int_{\R^{n}} |x|^{-b_a}U^{q-2}|v|^2~\!dx=\left(\frac{2q}{p}-1\right)I_{0}-2\beta H I_{1}.
\end{split}
\end{equation*}
Setting
$$
Q=\frac{q-p}{p-1}~\!,\quad K=\left(\frac{pH}{p-1}\right)^{p-2}C^{\frac{2q}{p}+p-2}\omega_{n-1}
$$
where $\omega_{n-1}=|\mathbb{S}^{n-1}|$ and $C$ is the constant in (\ref{eq:U_explicit}), and 
$$
\Phi(s,t)=K\int_{0}^{\infty}r^{s+n-1}\left(1+r^{QH}\right)^{-t}~\!dr\quad(0<s+n<tQH)
$$
we can write
\begin{equation*}
\begin{split}
&I_{0}=\left(\frac{pH}{p-1}\right)^{2}~\!\Phi\left(a+2\beta H-p+QHp,\frac{pq}{q-p}+2\right)\\
&I_{1}=\frac{pH}{p-1}~\!\Phi\left(a+2\beta H-p,\frac{pq}{q-p}+1\right)\\
&I_{2}=\Phi\left(a+2\beta H-p-QHp,\frac{pq}{q-p}\right)~\!.
\end{split}
\end{equation*}
As $0<s+n<tQH$, integration by parts yields that
$$
\Phi(s,t)=\frac{tQH}{s+n}~\!\Phi(s+QH,t+1).
$$
Hence
$$
I_{1}=\left(q-1+\frac{q}{p}\right)\frac{I_{0}}{s_{1}+n}\quad\text{and}\quad
I_{2}=\left(q-1+\frac{q}{p}\right)\frac{q~\!I_{0}}{(s_{1}+n)(s_{2}+n)}
$$
where
$$
s_{1}=a+2\beta H+(QH-1)p-QH\quad\text{and}\quad s_{2}=a+2\beta H+(QH-1)p-2QH~\!.
$$
Using the above computations and notation, from (\ref{eq:Q''}) we infer that
\begin{equation}
\label{eq:Q''beta}
\begin{split}
0&\le(s_{1}+n)(s_{2}+n)\left[\frac{(p-1)q^{2}}{p^{2}}-\left(\frac{2q}{p}-1\right)(q-1)\right]\\
&\quad-\frac{2(pq+q-p)}{p}\left[\frac{(p-1)q}{p}-(q-1)\right](s_{2}+n)\beta H\\
&\quad+\frac{q(pq+q-p)}{p}\left[(p-1)(\beta H)^{2}+(n-1)\right].
\end{split}
\end{equation}
The right hand side in (\ref{eq:Q''beta}) is a second order polynomial function with respect to $\beta$ having a minimum at
$$
\beta=\frac{Q}{p}=\frac{q-p}{p(p-1)}.
$$
We notice that this value of $\beta$ satisfies (\ref{eq:summability}). For such a value of $\beta$ (\ref{eq:Q''beta}) becomes
\begin{equation*}
\label{eq:BS2}
(n-1)H^{-2}\ge\frac{(q-p)(pq-q+p)}{p^{2}}~\!.
\end{equation*}
Hence if (\ref{eq:BS}) holds, then the extremal function for $S_{p}(a,q)$ cannot be radial.
\QED

\begin{Remark}
\label{R:v}
Let us focus our attention on the mapping $v=|x|^{\beta H}U^{\frac{q}{p}}$ used as a test function in the previous proof. When $p=2$ this choice of $v$ was already considered in \cite{FeSc}, 
 and in Proposition 8 of \cite{DELT},  to estimate the region of parameters $(a,q)$ for which symmetry breaking occurs. In \cite{FeSc} and in \cite{DELT} the expression for $v$ was suggested in a natural way by the explicit knowledge of the solutions of a related eigenvalue problem, according to known results displayed in the book \cite{LL} by Landau and Lifshitz. In the general case considered here, i.e., when $p\ne 2$, to our knowledge this information is missing. We finally point out that the exponent $\beta H$ in the expression of $v$ is the same appearing in the formula (\ref{eq:U_explicit}) of the radial ground state $U$. 
\end{Remark}
 
\section{Proof of Theorem \ref{T:2}}
\label{S:remarks}

We define
$$
S_{p}^{\rm rad}(a,q) :=\inf_{\scriptstyle u\in{\mathcal D}^{1,p}(\R^{n};|x|^a dx)\atop\scriptstyle u=u(|x|)~,~u\ne 0}\frac{\displaystyle\int_{\R^{n}}|x|^a|\nabla u|^p~\!dx}
{\left(\displaystyle\int_{\R^{n}} |x|^{-b_a}|u|^{q}~\!dx\right)^{p/q}}~\!.
$$
For a given pair $a,a'>p-n$ we put
$$
t=\frac{n-p+a'}{n-p+a}~\!.
$$
{From} the explicit knowledge of the extremals for $S_{p}^{\rm rad}(a,q)$
one infers that
\begin{equation}
\label{eq:rad}
S_{p}^{\rm rad}(a',q)= t^{p-1+\frac{p}{q}} S_{p}^{\rm rad}(a,q)~\!.
\end{equation}
Next we put
$$
\mathcal A:=\{a>p-n~|~S_{p}(a,q)<S_{p}^{\rm rad}(a,q)\}~\!.
$$
By Theorem \ref{T:BS1} the set $\mathcal A$ is not empty. 
Fix $a\in\mathcal A$ and take $a'>a$. 
We claim that $a'\in \mathcal A$. For the proof, let $u$ be an extremal for $S_{p}(a,q)$. Define $t$ as before, and put
$$
u_t(x)=u(|x|^{t-1}x)~\!.
$$
We compute, using polar coordinates $r=|x|, \sigma=\frac{x}{|x|}$,
\begin{eqnarray*}
\int_{\R^n}|x|^{a'}|\nabla u_t|^p~\!dx&=& t^{p-1}\int_{\R^n}|x|^a
\left|(\partial_r u)^2+t^{-2}|x|^{-2}|\nabla_\sigma u|^2\right|^{\frac{p}{2}}~\!dx\\
&\le&t^{p-1}\int_{\R^n}|x|^a
\left|(\partial_r u)^2+|x|^{-2}|\nabla_\sigma u|^2\right|^{\frac{p}{2}}~\!dx~\!,
\end{eqnarray*}
as $t>1$.
Since in addition
$$
\int_{\R^n}|x|^{-b_{a'}}|\nabla u_t|^q~\!dx= t^{-1}\int_{\R^n}|x|^{-b_a}
|u|^{q}~\!dx~\!,
$$
we infer that
\begin{eqnarray*}
S_{p}(a',q)&\le& \frac{\displaystyle\int_{\R^{n}}|x|^{a'}|\nabla u_t|^p~\!dx}
{\left(\displaystyle\int_{\R^{n}} |x|^{-b_{a'}}|u_t|^{q}~\!dx\right)^{p/q}}\le
 t^{p-1+\frac{p}{q}}\frac{\displaystyle\int_{\R^{n}}|x|^a|\nabla u|^p~\!dx}
{\left(\displaystyle\int_{\R^{n}} |x|^{-b_a}|u|^{q}~\!dx\right)^{p/q}}.
\end{eqnarray*}
Now we recall that $u\in\mathcal A$ achieves $S_{p}(a,q)$. Therefore, using also (\ref{eq:rad})
we conclude that
$$
S_{p}(a',q)\le t^{p-1+\frac{p}{q}}S_{p}(a,q)<t^{p-1+\frac{p}{q}}S_{p}^{\rm rad}(a,q)
=S_{p}^{\rm rad}(a',q),
$$
and therefore $a'\in\mathcal A$. We have proved that $\mathcal A\supseteq(\abs,\infty)$,
where $\abs:=\inf\mathcal A$. Actually,
$\mathcal A=(\abs,\infty)$, thanks to a continuity argument.

To check that
$a$ is nonnegative we first notice that
$S_{p}(a,q)\le S_{p}^{\rm rad}(a,q)$ for any $a>p-n$.
If $a=0$ the weight $|x|^{-b_0}$ is radially decreasing. Thus 
Schwarz symmetrization (see \cite{LiLo}) implies that 
$$
S_{p}(0,q)=\inf_{\scriptstyle u\in{\mathcal D}^{1,p}(\R^{n})
\atop\scriptstyle u\ne 0}\frac{\displaystyle\int_{\R^{n}}|\nabla u|^p~\!dx}{\left(\displaystyle\int_{\R^{n}} |x|^{-n+q\frac{n-p}{p}}|u|^{q}~\!dx\right)^{p/q}}=S_{p}^{\rm rad}(0,q)~\!,
$$
that is, $0\notin \mathcal A$ and therefore $0\le\abs$. 
\QED

\begin{Remark}
\label{R:np}
Assume $p\ge n$. The above argument gives the existence of $\abs\ge p-n$ such that symmetry breaking occurs in $(\abs,\infty)$ and can not occur in $(p-n,\abs]$. It would be of interest to know if $\abs>p-n$. A positive answer to this question was given by Dolbeault, Esteban, Loss and Tarantello in the Hilbertian case $n=p=2$, see Theorem 2 in \cite{DET}. 
\end{Remark}

\begin{Remark}
Here we assume that $1<p<n$ and we take $q=p^{*}$. From \cite{CKN} we know that the infimum
$$
S_{p}(a,p^{*}) :=\inf_{\scriptstyle u\in {\mathcal D}^{1,p}(\R^{n};|x|^a dx)\atop\scriptstyle u\ne 0}
\frac{\displaystyle\int_{\R^{n}}|x|^a|\nabla u|^p~\!dx}{\left(\displaystyle\int_{\R^{n}} |x|^{\frac{na}{n-p}}|u|^{p^{*}}~\!dx\right)^{p/p^{*}}}
$$
is positive, provided that $a>p-n$. Moreover, if in addition $a\le 0$ then $S_{p}(a,p^{*})$ is attained (see \cite{Au}, \cite{Ta} when $a=0$ and \cite{Ho}, \cite{GaMu1} when $a\in(p-n,0)$). 
The argument in the proof of Theorem \ref{T:2} has been already used by Horiuchi in \cite{Ho} to show that $S_{p}(a,p^{*})=S_{p}^{\rm rad}(a,p^{*})$ for any $a\le 0$.

If $a>0$ and $p=2$ then $S_{2}(a,2^{*})$ is not attained (see \cite{CaWa}) 
whereas if $p\ne 2$, to our knowledge any existence or non existence result of ground state is still missing. 
\end{Remark}

\newpage
\label{References}


\begin{thebibliography}{XX}
\footnotesize
%
%
\bibitem{Au} 
T. Aubin, Probl\`emes isop\'erim\'etriques et espaces de Sobolev, J. Differential Geometry {\bf 11} (1976), no.~4, 573--598. 

\bibitem{BW} 
J. Byeon\ and\ Z.-Q. Wang, Symmetry breaking of extremal functions for the Caffarelli-Kohn-Nirenberg inequalities, Commun. Contemp. Math. {\bf 4} (2002), no.~3, 457--465.

\bibitem{CKN}
L. Caffarelli, R. Kohn\ and\ L. Nirenberg, First order interpolation inequalities with weights, Compositio Math. {\bf 53} (1984), no.~3, 259--275.

\bibitem{CaMu99}
P. Caldiroli\ and\ R. Musina, On the existence of extremal functions for a weighted Sobolev embedding with critical exponent, Calc. Var. Partial Differential Equations {\bf 8} (1999), no.~4, 365--387.

\bibitem{CaWa}
F. Catrina\ and\ Z.-Q. Wang, On the Caffarelli-Kohn-Nirenberg inequalities: sharp constants, existence (and nonexistence), and symmetry of extremal functions, Comm. Pure Appl. Math. {\bf 54} (2001), no.~2, 229--258.

\bibitem{DEL}
J. Dolbeault, M. J. Esteban\ and\ M. Loss, Symmetry of extremals of functional inequalities via spectral estimates for linear operators, J. Math. Phys. {\bf 53} (2012), no.~9, 095204, 18 pp. 

\bibitem{DELT}
J. Dolbeault, M. J. Esteban, M. Loss\ and\  G. Tarantello, On the symmetry of extremals for the Caffarelli-Kohn-Nirenberg inequalities,
Adv. Nonlinear Stud. {\bf 9} (2009), no.~4, 713--726.

\bibitem{DET}
J. Dolbeault, M. J. Esteban\ and\ G. Tarantello, The role of Onofri type inequalities in the symmetry properties of extremals for Caffarelli-Kohn-Nirenberg inequalities, in two space dimensions, Ann. Sc. Norm. Super. Pisa Cl. Sci. (5) {\bf 7} (2008), no.~2, 313--341.

\bibitem{FeSc}
V. Felli\ and\ M. Schneider, Perturbation results of critical elliptic equations of Caffarelli-Kohn-Nirenberg type, J. Differential Equations {\bf 191} (2003), no.~1, 121--142.

\bibitem{GaMu1}
M. Gazzini\ and\ R. Musina, On a Sobolev-type inequality related to the weighted $p$-Laplace operator, J. Math. Anal. Appl. {\bf 352} (2009), no.~1, 99--111.

\bibitem{Ho}
T. Horiuchi, Best constant in weighted Sobolev inequality with weights being powers of distance from the origin, J. Inequal. Appl. {\bf 1} (1997), no.~3, 275--292.

\bibitem{LL}
{L. D. Landau\ and\ E. M. Lifshitz},
\textit{Quantum mechanics: non-relativistic theory. Theoretical Physics, Vol. 3}. 
Pergamon Press Ltd., London--Paris, 1958. 

\bibitem{LiLo}
E. H. Lieb\ and\ M. Loss, {\it Analysis}, second edition, Graduate Studies in Mathematics, 14, Amer. Math. Soc., Providence, RI, 2001. 

\bibitem{Maz}
V. G. Maz'ja, {\it Sobolev spaces}, translated from the Russian by T. O. Shaposhnikova, Springer Series in Soviet Mathematics, Springer, Berlin, 1985.

\bibitem{MuAMPA}
{R. Musina},
{Weighted Sobolev spaces of radially symmetric functions}, 
{Ann. Mat. Pura Appl.} (to appear) DOI: 10.1007/s10231-013-0348-4.

\bibitem{SW} 
D. Smets\ and\ M. Willem, Partial symmetry and asymptotic behavior for some elliptic variational problems, Calc. Var. Partial Differential Equations {\bf 18} (2003), no.~1, 57--75. 

\bibitem{Ta}
G. Talenti, Best constant in Sobolev inequality, Ann. Mat. Pura Appl. (4) {\bf 110} (1976), 353--372.
\end{thebibliography}
\end{document}